\documentclass[a4paper]{amsart}
\usepackage{color}
\usepackage[latin1]{inputenc}
\usepackage[T1]{fontenc}
\usepackage{amsfonts}
\usepackage{amssymb}
\usepackage{amsmath}
\usepackage{amsthm}
\usepackage{stmaryrd}
\usepackage{enumerate}
\usepackage{multirow}
\usepackage{graphicx}
\usepackage{setspace}

\usepackage{graphicx,type1cm,eso-pic,color}
\usepackage{pstricks}
\usepackage{float}
\newtheorem{theorem}{Theorem}[section]
\newtheorem{lemma}[theorem]{Lemma}

\begin{document}
\setlength\arraycolsep{2pt}
\title{Multifunction Estimation in a Time-Discretized Skorokhod Reflection Problem}
\author{Nicolas MARIE$^{\dag}$}
\email{nmarie@parisnanterre.fr}
\date{}
\maketitle
\noindent
$^{\dag}$Universit\'e Paris Nanterre, CNRS, Modal'X, 92001 Nanterre, France.
%


%
\begin{abstract}
This paper deals with a consistent estimator of the multifunction involved in a time-discretized Skorokhod reflection problem defined by a stochastic differential equation and a Moreau sweeping process.
\end{abstract}
%


%
\section{Introduction}\label{section_introduction}
Let $\texttt X_t$ be the location recorded at time $t\in [0,T]$ of a target living in a convex body $\texttt C(t)$ of $\mathbb R^m$, where $T > 0$ and $m\in\mathbb N^*$. Due to a lack of information on the target's behavior, $\texttt C(t)$ may be unknown, and then needs to be estimated. For instance, if $\texttt X_t$ is the location of an animal at time $t$, then its territory $\texttt C(t)$ is a priori unknown. Such an application in ethology has been investigated in Cholaquidis et al. \cite{CFLP16} when the process $\texttt X$ is a reflected Brownian motion. More generally, $\texttt X$ may be modeled by the following continuous-time Skorokhod reflection problem defined by a stochastic differential equation and a Moreau sweeping process (see Moreau \cite{MOREAU76}):
\begin{eqnarray}
 \texttt X_t & = &
 \int_{0}^{t}b(\texttt X_s)ds +
 \int_{0}^{t}\sigma(\texttt X_s)d\texttt W_s +\texttt Y_t
 \nonumber\\
 \label{Skorokhod_reflection_problem}
 & &
 \hspace{2cm}{\rm with}\quad
 \left\{
 \begin{array}{rcl}
  -\dot{\tt Y}_t & \in & \mathcal N_{\texttt C(t)}(\texttt X_t)\quad
  |\textrm D\texttt Y|\textrm{-a.e.}\\
  \texttt Y_0(\cdot) & \equiv & \mathbf x_0
 \end{array}
 \right.\textrm{$;$ }t\in [0,T],
\end{eqnarray}
where
\begin{itemize}
 \item $b$ (resp. $\sigma$) is a Lipschitz continuous map from $\mathbb R^m$ into itself (resp. ${\rm GL}_m(\mathbb R)$; the space of the $m\times m$ invertible matrices with real entries), $\texttt W$ is a $m$-dimensional Brownian motion and $\texttt C :\mathbb R_+\rightrightarrows\mathbb R^m$ is a convex-body-valued multifunction.
 \item For every convex body $K$ of $\mathbb R^m$,
 \begin{displaymath}
 \mathcal N_K(\mathbf x) :=
 \{\mathbf s\in\mathbb R^m :\forall\mathbf y\in K\textrm{, }
 \langle\mathbf s,\mathbf y -\mathbf x\rangle\leqslant 0\}
 \quad\textrm{is the normal cone of $K$ at point $\mathbf x\in\mathbb R^m$.}
 \end{displaymath}
 \item For every continuous function of bounded variation $\varphi : [0,T]\rightarrow\mathbb R^m$,
 \begin{displaymath}
 \dot\varphi :=\frac{d\textrm D\varphi}{d|\textrm D\varphi|}
 \quad\textrm{
 with $\textrm D\varphi$ the differential measure of $\varphi$ and $|\textrm D\varphi|$ its variation measure.}
 \end{displaymath}
\end{itemize}
Roughly speaking, a solution of Problem (\ref{Skorokhod_reflection_problem}) behaves as a diffusion process when $\texttt X_t\in {\rm int}(\texttt C(t))$, but when $\texttt X_t$ reaches $\partial\texttt C(t)$, a minimal force $\texttt Y_t$ is applied to $\texttt X_t$ in order to keep it inside of $\texttt C(t)$. About the existence and uniqueness of the solution of Problem (\ref{Skorokhod_reflection_problem}), the reader may refer to Bernicot and Venel \cite{BV11} and Castaing et al. \cite{CMR16}. Under appropriate regularity conditions on the multifunction $\texttt C$ (see Bernicot and Venel \cite{BV11}, Theorems 3.9 and 4.6), a uniformly convergent approximation of $\texttt X$ is defined by
\begin{equation}\label{Euler_scheme}
\left\{
\begin{array}{rcl}
 \overline{\tt X}_{j + 1} & = & \Pi_{\texttt C(t_{j + 1})}(
 \overline{\tt X}_j +
 b(\overline{\tt X}_j)\Delta_n +
 \sigma(\overline{\tt X}_j)(\texttt W_{t_{j + 1}} -\texttt W_{t_j}))\\
 \overline{\tt X}_0 & \equiv & \mathbf x_0\in\texttt C(0)
\end{array}\right.
\textrm{$;$ }j\in J_n\backslash\{n\},
\end{equation}
where $\Pi_C(\cdot)$ is the orthogonal projection from $\mathbb R^m$ onto a given convex body $C$ of $\mathbb R^m$, $J_n :=\{0,\dots,n\}$ with $n\in\mathbb N^*$, $\Delta_n := T/n$ with $T > 0$, and $t_j := j\Delta_n$ for every $j\in J_n$. Our paper, which is part of the copies-based statistical inference for diffusion processes (see Comte and Genon-Catalot \cite{CGC20}, Marie and Rosier \cite{MR23}, Denis et al. \cite{DDM21}, etc.), deals with the consistency of the estimator
\begin{displaymath}
\widehat{\tt C}_{N,j} :=
{\rm conv}\{\overline{\tt X}_{j}^{1},\dots,\overline{\tt X}_{j}^{N}\}
\textrm{ $;$ }j\in J_n\backslash\{0\}
\end{displaymath}
of $\texttt C(t_j)$, where ${\rm conv}(\cdot)$ is the convex hull operator, and $\overline{\tt X}^1,\dots,\overline{\tt X}^N$ are $N$ independent copies of $\overline{\tt X}$. Precisely, Theorem \ref{convergence_unidimensional} says that $\widehat{\tt C}_{N,j}$ converges to $\texttt C(t_j)$ - with rate $N$ - when $m = 1$, and Theorem \ref{convergence_multidimensional} says that $\widehat{\tt C}_{N,j}$ converges pointwise to $\texttt C(t_j)$ when $m > 1$ and the multifunction $\texttt C$ is decreasing (possibly constant). Finally, for another example of estimator computed from observations of the approximation scheme of a reflected diffusion model, the reader may refer to Cattiaux et al. \cite{CLP17}.
\\
\\
{\bf Notations and basic definitions:}
\begin{itemize}
 \item The Euclidean norm (resp. distance) on $\mathbb R^m$ is denoted by $\|.\|$ (resp. $d$).
 \item For every $m\times m$ real matrix $\mathbf A$,
 \begin{displaymath}
 \|\mathbf A\|_{\rm op} :=
 \sup_{\mathbf x\in\mathbb R^m}\frac{\|\mathbf A\mathbf x\|}{\|\mathbf x\|}.
 \end{displaymath}
 \item For every $K\in\mathcal P(\mathbb R^m)$,
 \begin{displaymath}
 d(\cdot,K) :=
 \inf\{d(\cdot,\mathbf y)\textrm{ $;$ }\mathbf y\in K\}.
 \end{displaymath}
 \item Consider two normed vector spaces $(E,\|.\|_E)$ and $(F,\|.\|_F)$. For every Lipschitz continuous map $f : E\rightarrow F$,
 \begin{displaymath}
 \|f\|_{\rm Lip} :=
 \sup_{x\neq y}\frac{\|f(y) - f(x)\|_F}{\|y - x\|_E}.
 \end{displaymath}
\end{itemize}
%


%
\section{Main results}\label{section_main_results}
First, the following theorem says that $\widehat{\tt C}_{N,j}$ converges to $\texttt C(t_j)$ - with rate $N$ - when $m = 1$.
%


%
\begin{theorem}\label{convergence_unidimensional}
Assume that there exists a constant $\mathfrak m_{\tt C} > 0$ such that, for every $t\in [0,T]$ and $x\in {\tt C}(t)$, $|x|\leqslant\mathfrak m_{\tt C}$. For any $j\in J_n\backslash\{0\}$,
\begin{displaymath}
d_{\rm H}(\widehat{\tt C}_{N,j},{\tt C}(t_j)) :=
\max\{\sup({\tt C}(t_j)) -\sup(\widehat{\tt C}_{N,j}),
\inf(\widehat{\tt C}_{N,j}) -\inf({\tt C}(t_j))\}
\xrightarrow[N\rightarrow\infty]{\mathbb P} 0,
\end{displaymath}
and if in addition $\sigma$ is constant, then
\begin{displaymath}
Nd_{\rm H}(\widehat{\tt C}_{N,j},{\tt C}(t_j))\xrightarrow[N\rightarrow\infty]{\mathbb P} 0.
\end{displaymath}
\end{theorem}
\noindent
Now, consider two convex bodies $K$ and $L$ of $\mathbb R^m$ such that $L\subset K$. If
\begin{displaymath}
d(\mathbf x,L) = 0
\textrm{ $;$ }
\forall\mathbf x\in {\rm int}(K),
\end{displaymath}
then ${\rm int}(K)\subset L$, leading to $K = L$ because $K$ and $L$ are closed subsets of $\mathbb R^m$. For this reason, it makes sense to say that a sequence $(L_n)_{n\in\mathbb N}$ of random convex-compact subsets of $K$ converges pointwise (in probability) to $K$ if and only if,
\begin{displaymath}
d(\mathbf x,L_n)\xrightarrow[n\rightarrow\infty]{\mathbb P} 0
\textrm{ $;$ }\forall\mathbf x\in {\rm int}(K).
\end{displaymath}
The following theorem says that $\widehat{\tt C}_{N,j}$ is pointwise consistent, for every $j\in J_n\backslash\{0\}$, when the multifunction $\texttt C$ is decreasing (possibly constant).
%


%
\begin{theorem}\label{convergence_multidimensional}
Assume that ${\tt C}(t)\subset {\tt C}(s)$ for every $s,t\in [0,T]$ satisfying $t\geqslant s$. Then, for every $j\in J_n\backslash\{0\}$, $\widehat{\tt C}_{N,j}$ converges pointwise to ${\tt C}(t_j)$:
\begin{displaymath}
d(\mathbf x,\widehat{\tt C}_{N,j})
\xrightarrow[N\rightarrow\infty]{\mathbb P} 0
\textrm{ $;$ }\forall\mathbf x\in {\rm int}({\tt C}(t_j)).
\end{displaymath}
\end{theorem}
\appendix
%


%
\section{Proofs}\label{section_proofs}
%


%
\subsection{Proof of Theorem \ref{convergence_unidimensional}}
The proof of Theorem \ref{convergence_unidimensional} relies on the following lemma.
%


%
\begin{lemma}\label{lemma_convergence_unidimensional}
Consider $C = [I,S]$ with $I,S\in\mathbb R$ satisfying $I < S$, a $\mathbb R$-valued random variable $H$, and
\begin{eqnarray*}
 \widehat C_N & := & [\widehat I_N,\widehat S_N]
 \quad {\rm with}\quad
 \widehat I_N =
 \min\{X_1,\dots,X_N\}\\
 & &
 \hspace{5cm}{\rm and}\quad
 \widehat S_N =
 \max\{X_1,\dots,X_N\},
\end{eqnarray*}
where $X_1,\dots,X_N$ are $N$ independent copies of $X =\Pi_C(H)$. Moreover, let $\Phi$ (resp. $F$) be the distribution function of $H$ (resp. $X$).
\begin{enumerate}
 \item Assume that
 \begin{equation}\label{lemma_convergence_unidimensional_1}
 \Phi((I,S))\subset (0,1).
 \end{equation}
 Then,
 \begin{displaymath}
 d_{\rm H}(\widehat C_N,C)\xrightarrow[N\rightarrow\infty]{\mathbb P} 0.
 \end{displaymath}
 \item Assume that there exist $\varphi,\psi\in C^1(C;(0,1))$ such that
 \begin{equation}\label{lemma_convergence_unidimensional_2}
 \psi(\cdot)\leqslant\Phi(\cdot)\leqslant\varphi(\cdot)
 \quad\textrm{on}\quad C.
 \end{equation}
 Then,
 \begin{displaymath}
 Nd_{\rm H}(\widehat C_N,C)\xrightarrow[N\rightarrow\infty]{\mathbb P} 0.
 \end{displaymath}
\end{enumerate}
\end{lemma}
\noindent
The proof of Lemma \ref{lemma_convergence_unidimensional} is postponed to Subsubsection \ref{section_proof_lemma_convergence_unidimensional}.
\\
\\
Consider $j\in J_n\backslash\{n\}$, and let $\Phi_{j + 1}$ be the distribution function of
\begin{displaymath}
\texttt H_{j + 1} :=
\overline{\tt X}_j + b(\overline{\tt X}_j)\Delta_n +
\sigma(\overline{\tt X}_j)\texttt Z_{j + 1}
\quad {\rm with}\quad
\texttt Z_{j + 1} =\texttt W_{t_{j + 1}} -\texttt W_{t_j}.
\end{displaymath}
First, since $\overline{\tt X}_j\in\texttt C(t_j)$, $|\overline{\tt X}_j|\leqslant\mathfrak m_{\tt C}$, and then
\begin{displaymath}
\underbrace{-\mathfrak m_{\tt C} +
\Delta_n\min_{|v|\leqslant\mathfrak m_{\tt C}}b(v)}_{=: m_n} +
\sigma(\overline{\tt X}_j)\texttt Z_{j + 1}
\leqslant\texttt H_{j + 1}\leqslant
\underbrace{\mathfrak m_{\tt C} +
\Delta_n\max_{|v|\leqslant\mathfrak m_{\tt C}}b(v)}_{=: M_n} +
\sigma(\overline{\tt X}_j)\texttt Z_{j + 1}.
\end{displaymath}
So, by assuming that $\sigma(\cdot) > 0$ without loss of generality, for any $x\in\mathbb R$,
\begin{eqnarray}
 \mathbb P\left(\texttt Z_{j + 1}\leqslant -\frac{|x - M_n|}{m_{\sigma}}\right)
 & \leqslant &
 \mathbb P(\sigma(\overline{\tt X}_j)\texttt Z_{j + 1}\leqslant x - M_n)
 \nonumber\\
 & \leqslant &
 \Phi_{j + 1}(x)
 \nonumber\\
 \label{convergence_unidimensional_1}
 & \leqslant &
 \mathbb P(\sigma(\overline{\tt X}_j)\texttt Z_{j + 1}\leqslant x - m_n)
 \leqslant
 \mathbb P\left(\texttt Z_{j + 1}\leqslant\frac{|x - m_n|}{m_{\sigma}}\right)
\end{eqnarray}
with
\begin{displaymath}
m_{\sigma} =\min_{|x|\leqslant\mathfrak m_{\tt C}}\sigma(x).
\end{displaymath}
Since $\texttt Z_{j + 1}$ is a Gaussian random variable, $\Phi_{j + 1}(x)\in (0,1)$, and then
\begin{displaymath}
d_{\rm H}(\widehat{\tt C}_{N,j},{\tt C}(t_j))\xrightarrow[N\rightarrow\infty]{\mathbb P} 0
\quad\textrm{by Lemma \ref{lemma_convergence_unidimensional}.(1).}
\end{displaymath}
Now, assume that $\sigma(\cdot)\equiv\texttt s$ with $\texttt s > 0$. By Inequality (\ref{convergence_unidimensional_1}),
\begin{displaymath}
\psi_{j + 1}(\cdot)\leqslant
\Phi_{j + 1}(\cdot)\leqslant
\varphi_{j + 1}(\cdot)
\quad {\rm on}\quad\mathbb R,
\end{displaymath}
where $\varphi_{j + 1}$ (resp. $\psi_{j + 1}$) is the distribution function of $\xi_{j + 1} :=\texttt s\texttt Z_{j + 1} + M_n$ (resp. $\zeta_{j + 1} :=\texttt s\texttt Z_{j + 1} + m_n$). Since both $\xi_{j + 1}$ and $\zeta_{j + 1}$ are Gaussian random variables, $\varphi_{j + 1},\psi_{j + 1}\in C^1(\mathbb R;(0,1))$, and then
\begin{displaymath}
Nd_{\rm H}(\widehat{\tt C}_{N,j},{\tt C}(t_j))\xrightarrow[N\rightarrow\infty]{\mathbb P} 0
\quad\textrm{by Lemma \ref{lemma_convergence_unidimensional}.(2).}
\end{displaymath}
%


%
\subsubsection{Proof of Lemma \ref{lemma_convergence_unidimensional}}\label{section_proof_lemma_convergence_unidimensional}
First of all, in order to prove Lemma \ref{lemma_convergence_unidimensional}, let us establish the following suitable relationship between $\Phi$ and $F$:
\begin{equation}\label{lemma_convergence_unidimensional_3}
F(x) =\Phi(x)\mathbf 1_{I\leqslant x < S} +\mathbf 1_{x\geqslant S}
\textrm{ $;$ }\forall x\in\mathbb R.
\end{equation}
For any $x\in\mathbb R$, since $\{I < H < S\}\subset \{X = H\}$,
\begin{eqnarray*}
 F(x) & = &
 \mathbb P(X\leqslant x,H\in (I,S)) +
 \mathbb P(X\leqslant x,H\not\in (I,S))\\
 & = &
 \mathbb P(H\leqslant x,I < H < S) +
 \mathbb P(X\leqslant x,H\leqslant I) +
 \mathbb P(X\leqslant x,H\geqslant S).
\end{eqnarray*}
Moreover,
\begin{displaymath}
\mathbb P(H\leqslant x,I < H < S) =
\mathbb P(I < H\leqslant x)\mathbf 1_{I < x < S} +
\mathbb P(I < H < S)\mathbf 1_{x\geqslant S}
\end{displaymath}
and, since $\{H\leqslant I\} =\{X = I\}$ and $\{H\geqslant S\} =\{X = S\}$,
\begin{displaymath}
\mathbb P(X\leqslant x,H\leqslant I) +
\mathbb P(X\leqslant x,H\geqslant S) =
\mathbb P(H\leqslant I)\mathbf 1_{x\geqslant I} +
\mathbb P(H\geqslant S)\mathbf 1_{x\geqslant S}.
\end{displaymath}
Thus,
\begin{eqnarray*}
 F(x) & = &
 \mathbb P(I < H\leqslant x)\mathbf 1_{I < x < S} +
 \mathbb P(H\leqslant I)\mathbf 1_{x\geqslant I}\\
 & &
 \hspace{2.5cm} +
 (\mathbb P(I < H < S) +
 \mathbb P(H\geqslant S))\mathbf 1_{x\geqslant S}\\
 & = &
 (\mathbb P(H\leqslant x) -\mathbb P(H\leqslant I))\mathbf 1_{I < x < S}\\
 & &
 \hspace{2.5cm} +
 \mathbb P(H\leqslant I)(\mathbf 1_{x = I} +\mathbf 1_{x > I}) +
 \mathbb P(H > I)\mathbf 1_{x\geqslant S}\\
 & = &
 \mathbb P(H\leqslant x)\mathbf 1_{I\leqslant x < S} +
 \underbrace{(\mathbb P(H\leqslant I) +\mathbb P(H > I))}_{= 1}
 \mathbf 1_{x\geqslant S} =
 \Phi(x)\mathbf 1_{I\leqslant x < S} +\mathbf 1_{x\geqslant S}.
\end{eqnarray*}
Let us establish Lemma \ref{lemma_convergence_unidimensional}.(1,2).
\begin{enumerate}
 \item First, let us show that $H$ satisfies (\ref{lemma_convergence_unidimensional_1}) if and only if
 \begin{equation}\label{lemma_convergence_unidimensional_4}
 \Phi((-\infty,S))\subset [0,1)\quad {\rm and}\quad
 \Phi((I,\infty))\subset (0,1].
 \end{equation}
 On the one hand, if $H$ satisfies (\ref{lemma_convergence_unidimensional_4}), then
 \begin{eqnarray*}
  \Phi((I,S)) & = &
  \Phi((-\infty,S)\cap (I,\infty))\\
  & \subset &
  \Phi((-\infty,S))\cap\Phi((I,\infty))
  \subset [0,1)\cap (0,1] = (0,1).
 \end{eqnarray*}
 On the other hand, (\ref{lemma_convergence_unidimensional_1}) implies (\ref{lemma_convergence_unidimensional_4}) because $\Phi$ is increasing. Now, let us prove that if
 \begin{equation}\label{lemma_convergence_unidimensional_5}
 \Phi((-\infty,S))\subset [0,1),
 \end{equation}
 then $\widehat S_N$ is a converging estimator of $S$. By (\ref{lemma_convergence_unidimensional_3}) and (\ref{lemma_convergence_unidimensional_5}),
 \begin{equation}\label{lemma_convergence_unidimensional_6}
 \sup\{x\in\mathbb R : F(x) < 1\} = S.
 \end{equation}
 Let $F_N$ be the distribution function of $\widehat S_N$. Since $X_1,\dots,X_N$ are independent copies of $X$, for any $x\in\mathbb R$,
 \begin{displaymath}
 F_N(x) =
 \mathbb P\left(\bigcap_{i = 1}^{N}\{X_i\leqslant x\}\right) = F(x)^N.
 \end{displaymath}
 On the one hand, for $x < S$, $F(x)\in [0,1)$ by (\ref{lemma_convergence_unidimensional_6}), and then
 \begin{displaymath}
 \lim_{N\rightarrow\infty}F_N(x) =
 \lim_{N\rightarrow\infty}F(x)^N = 0.
 \end{displaymath}
 On the other hand, for $x\geqslant S$, $F(x) = 1$ by (\ref{lemma_convergence_unidimensional_6}), and then
 \begin{displaymath}
 \lim_{N\rightarrow\infty}F_N(x) =
 \lim_{N\rightarrow\infty}F(x)^N = 1.
 \end{displaymath}
 So, $F_N$ converges pointwise to $\mathbf 1_{[S,\infty)}$ (the distribution function of $\omega\mapsto S$), leading to
 \begin{displaymath}
 \widehat S_N\xrightarrow[N\rightarrow\infty]{\mathbb P}S.
 \end{displaymath}
 Finally, by assuming that $H$ satisfies (\ref{lemma_convergence_unidimensional_1}) (or equivalently (\ref{lemma_convergence_unidimensional_4})),
 \begin{displaymath}
 \widehat I_N\xrightarrow[N\rightarrow\infty]{\mathbb P}I,
 \end{displaymath}
 and thus
 \begin{eqnarray*}
  d_{\rm H}(\widehat C_N,C) & = &
  \sup_{x\in C}d(x,\widehat C_N)
  \quad {\rm because}\quad
  \widehat C_N\subset C\\
  & = &
  \max\{S -\widehat S_N,\widehat I_N - I\}
  \xrightarrow[N\rightarrow\infty]{\mathbb P} 0.
 \end{eqnarray*}
 \item First, since $X$ is a $[I,S]$-valued random variable, for every $x\in\mathbb R_+$,
 \begin{displaymath}
 \mathbb P(N(\widehat S_N - S)\leqslant x) =
 \mathbb P\left(X\leqslant S +\frac{x}{N}\right)^N = 1.
 \end{displaymath}
 Now, by Equality (\ref{lemma_convergence_unidimensional_3}), and since $\Phi(\cdot)\leqslant\varphi(\cdot)$ on $C = [I,S]$ with $\varphi\in C^1(C;(0,1))$, for every $y\in [I,S)$,
 \begin{eqnarray*}
  1 -\mathbb P(X > y) = F(y)
  & = &
  \Phi(y)\\
  & \leqslant &
  \varphi(y)
  \underset{y\uparrow S}{=}
  \varphi(S) +\varphi'(S)(y - S) +\textrm o(y - S).
 \end{eqnarray*}
 Then, for any $x < 0$, and every $N\in\mathbb N^*$ such that $N\geqslant -x/(S - I)$,
 \begin{eqnarray*}
  \mathbb P(N(\widehat S_N - S)\leqslant x)
  & = &
  \left(1 -\mathbb P\left(X > S +\frac{x}{N}\right)\right)^N\\
  & \leqslant &
  \varphi\left(S +\frac{x}{N}\right)^N =
  \varphi(S)^N\exp(N\log(1 + U_N(x)))\\
  & &
  \hspace{4.5cm}{\rm with}\quad
  U_N(x) =\frac{\varphi(S + x/N)}{\varphi(S)} - 1,
 \end{eqnarray*}
 and
 \begin{displaymath}
 N\log(1 + U_N(x))
 \underset{N\rightarrow\infty}{=}
 \frac{\varphi'(S)}{\varphi(S)}x +\textrm o(1).
 \end{displaymath}
 Since $\varphi(S) < 1$, the distribution function of $N(\widehat S_N - S)$ converges pointwise to $\mathbf 1_{\mathbb R_+}$, and thus
 \begin{equation}\label{lemma_convergence_unidimensional_7}
 N(S -\widehat S_N)
 \xrightarrow[N\rightarrow\infty]{\mathbb P} 0.
 \end{equation}
 Since $\Phi(\cdot)\geqslant\psi(\cdot)$ on $C = [I,S]$ with $\psi\in C^1(C;(0,1))$, as previously,
 \begin{equation}\label{lemma_convergence_unidimensional_8}
 N(\widehat I_N - I)
 \xrightarrow[N\rightarrow\infty]{\mathbb P} 0.
 \end{equation}
 In conclusion, by (\ref{lemma_convergence_unidimensional_7}) and (\ref{lemma_convergence_unidimensional_8}),
 \begin{displaymath}
 Nd_{\rm H}(\widehat C_N,C) =
 \max\{N(S -\widehat S_N),N(\widehat I_N - I)\}
 \xrightarrow[N\rightarrow\infty]{\mathbb P} 0.
 \end{displaymath}
\end{enumerate}
%


%
\subsection{Proof of Theorem \ref{convergence_multidimensional}}\label{section_proof_convergence_multidimensional}
The proof of Theorem \ref{convergence_multidimensional} relies on the following lemma.
%


%
\begin{lemma}\label{lemma_convergence_multidimensional}
Consider a convex body $C$ of $\mathbb R^m$, a $\mathbb R^m$-valued random variable $H$, and
\begin{displaymath}
\widehat C_N := {\rm conv}\{X_1,\dots,X_N\},
\end{displaymath}
where $X_1,\dots,X_N$ are $N$ independent copies of $X =\Pi_C(H)$. Assume that, for every $\mathbf x\in {\rm int}(C)$ and $\varepsilon > 0$,
\begin{equation}\label{lemma_convergence_multidimensional_1}
\mathbb P_H(\overline B_d(\mathbf x,\varepsilon)\cap {\rm int}(C)) > 0.
\end{equation}
Then,
\begin{displaymath}
d(\mathbf x,\widehat C_N)\xrightarrow[N\rightarrow\infty]{\mathbb P} 0
\textrm{ $;$ }\forall\mathbf x\in {\rm int}(C).
\end{displaymath}
\end{lemma}
\noindent
The proof of Lemma \ref{lemma_convergence_multidimensional} is postponed to Subsubsection \ref{section_proof_lemma_convergence_multidimensional}.
\\
\\
For every $j\in J_n\backslash\{n\}$, consider
\begin{displaymath}
\texttt H_{j + 1} :=
\overline{\tt X}_j + b(\overline{\tt X}_j)\Delta_n +
\sigma(\overline{\tt X}_j)\texttt Z_{j + 1}
\quad {\rm with}\quad
\texttt Z_{j + 1} =\texttt W_{t_{j + 1}} -\texttt W_{t_j}.
\end{displaymath}
The proof of Theorem \ref{convergence_multidimensional} is dissected in three steps. Step 1 deals with a suitable control of $\|\texttt H_{j + 1} -\mathbf x\|$, depending on $\texttt Z_{j + 1}$ and $\|\overline{\tt X}_j -\mathbf x\|$, for every $j\in J_n\backslash\{n\}$ and $\mathbf x\in\texttt C(t_{j + 1})$. In Step 2, the condition
\begin{eqnarray*}
 & &
 \mathbb P_{\texttt H_j}(\overline B_d(\mathbf x,\varepsilon)
 \cap {\rm int}(\texttt C(t_j))) > 0\textrm{ $;$}\\
 & &
 \hspace{3cm}
 \forall j\in J_n\backslash\{0\}\textrm{, }
 \forall\varepsilon > 0\textrm{, }
 \forall\mathbf x\in {\rm int}(\texttt C(t_j))
\end{eqnarray*}
is checked in order to apply Lemma \ref{lemma_convergence_multidimensional} to $\widehat{\tt C}_{N,1},\dots,\widehat{\tt C}_{N,n}$ in Step 3.
\\
\\
{\bf Step 1.} For any $j\in J_n\backslash\{n\}$ and $\mathbf x\in\texttt C(t_{j + 1})$, since $b$ and $\sigma$ are Lipschitz continuous maps,
\begin{eqnarray*}
 \|\texttt H_{j + 1} -\mathbf x\| & = &
 \|\overline{\tt X}_j -\mathbf x + (b(\overline{\tt X}_j) - b(\mathbf x))\Delta_n\\
 & &
 \hspace{2cm} +
 (\sigma(\overline{\tt X}_j) -\sigma(\mathbf x))\texttt Z_{j + 1} +
 b(\mathbf x)\Delta_n +\sigma(\mathbf x)\texttt Z_{j + 1}\|\\
 & \leqslant &
 (1 +\|b\|_{\rm Lip}\Delta_n +\|\sigma\|_{\rm Lip}\|\texttt Z_{j + 1}\|)
 \|\overline{\tt X}_j -\mathbf x\| +\texttt R_{j + 1}(\mathbf x)\\
 & &
 \hspace{4cm}{\rm with}\quad
 \texttt R_{j + 1}(\mathbf x) =
 \|\sigma(\mathbf x)\texttt Z_{j + 1} + b(\mathbf x)\Delta_n\|.
\end{eqnarray*}
On the one hand, since $\sigma(\mathbf x)\in {\rm GL}_m(\mathbb R)$,
\begin{eqnarray*}
 \|\texttt Z_{j + 1}\| & = &
 \|\sigma(\mathbf x)^{-1}(\sigma(\mathbf x)\texttt Z_{j + 1} +
 b(\mathbf x)\Delta_n -b(\mathbf x)\Delta_n)\|\\
 & \leqslant &
 \|\sigma(\mathbf x)^{-1}\|_{\rm op}(
 \texttt R_{j + 1}(\mathbf x) +\|b(\mathbf x)\|\Delta_n).
\end{eqnarray*}
On the other hand, since $\texttt C(t_0),\dots,\texttt C(t_n)\subset\texttt C(0)$, and since $\texttt C(0)$ is a compact subset of $\mathbb R^m$, there exists a constant $\mathfrak m_{\tt C} > 0$ such that, for every $k\in J_n$ and $\mathbf v\in\texttt C(t_k)$, $\|\mathbf v\|\leqslant\mathfrak m_{\tt C}$. Thus,
\begin{equation}\label{convergence_multidimensional_1}
\|\texttt H_{j + 1} -\mathbf x\|
\leqslant
\mathfrak c_1\|\overline{\tt X}_j -\mathbf x\| +
\mathfrak c_2\texttt R_{j + 1}(\mathbf x)
\end{equation}
with
\begin{eqnarray*}
 \mathfrak c_1 & = &
 1 +\left(\|b\|_{\rm Lip} +
 \|\sigma\|_{\rm Lip}
 \sup_{\|\mathbf v\|\leqslant\mathfrak m_{\tt C}}\{
 \|\sigma(\mathbf v)^{-1}\|_{\rm op}\|b(\mathbf v)\|\}\right)\Delta_n\\
 & &
 \hspace{4cm}{\rm and}\quad
 \mathfrak c_2 = 1 +
 2\mathfrak m_{\tt C}\|\sigma\|_{\rm Lip}
 \sup_{\|\mathbf v\|\leqslant\mathfrak m_{\tt C}}
 \|\sigma(\mathbf v)^{-1}\|_{\rm op}.
\end{eqnarray*}
{\bf Step 2.} For every $j\in J_n\backslash\{0\}$ and $\mathbf x\in {\rm int}(\texttt C(t_j))$, there exists $\varepsilon_{j,\mathbf x} > 0$ such that
\begin{equation}\label{convergence_multidimensional_2}
\overline B_d(\mathbf x,\varepsilon_{j,\mathbf x})\subset
{\rm int}(\texttt C(t_j)).
\end{equation}
First, for every $\varepsilon > 0$ and $\mathbf x\in {\rm int}(\texttt C(t_1))$, by (\ref{convergence_multidimensional_2}), since $\sigma(\mathbf x_0)\in {\rm GL}_m(\mathbb R)$, and since $[\texttt W_{t_1}]_1,\dots,[\texttt W_{t_1}]_m$ are i.i.d. Gaussian random variables,
\begin{eqnarray*}
 & &
 \mathbb P_{\texttt H_1}(\overline B_d(\mathbf x,\varepsilon)
 \cap {\rm int}(\texttt C(t_1)))\\
 & &
 \hspace{1cm}\geqslant
 \mathbb P(\|\mathbf x_0 -\mathbf x +
 b(\mathbf x_0)\Delta_n +
 \sigma(\mathbf x_0)\texttt W_{t_1}\|
 \leqslant\varepsilon\wedge\varepsilon_{1,\mathbf x})\\
 & &
 \hspace{1cm}\geqslant
 \mathbb P(\|\texttt W_{t_1} +
 \sigma(\mathbf x_0)^{-1}(\mathbf x_0 -\mathbf x +
 b(\mathbf x_0)\Delta_n)\|
 \leqslant\|\sigma(\mathbf x_0)\|_{\rm op}^{-1}
 (\varepsilon\wedge\varepsilon_{1,\mathbf x}))
 \geqslant\mu_{1,\mathbf x}(\varepsilon)
\end{eqnarray*}
with
\begin{eqnarray*}
 \mu_{1,\mathbf x}(\varepsilon) & = &
 \prod_{k = 1}^{m}
 \mathbb P(|[\texttt W_{t_1}]_k +
 [\sigma(\mathbf x_0)^{-1}(\mathbf x_0 -\mathbf x +
 b(\mathbf x_0)\Delta_n)]_k|\\
 & &
 \hspace{5cm}\leqslant
 (\sqrt m\|\sigma(\mathbf x_0)\|_{\rm op})^{-1}
 (\varepsilon\wedge\varepsilon_{1,\mathbf x})) > 0.
\end{eqnarray*}
Now, consider $j\in J_n\backslash\{0,n\}$, and assume that for every $\varepsilon > 0$ and $\mathbf x\in {\rm int}(\texttt C(t_j))$,
\begin{equation}\label{convergence_multidimensional_3}
\mathbb P_{\texttt H_j}(\overline B_d(\mathbf x,\varepsilon)
\cap {\rm int}(\texttt C(t_j))) > 0.
\end{equation}
For any $\varepsilon > 0$ and $\mathbf x\in {\rm int}(\texttt C(t_{j + 1}))$, by Inequality (\ref{convergence_multidimensional_1}), by (\ref{convergence_multidimensional_2}), since $\texttt Z_{j + 1}$ (and then $\texttt R_{j + 1}$) is independent of $\mathcal F_j =\sigma(\texttt W_{t_0},\dots,\texttt W_{t_j})$, and since $\overline{\tt X}_j$ is $\mathcal F_j$-measurable,
\begin{eqnarray*}
 & &
 \mathbb P_{\texttt H_{j + 1}}(
 \overline B_d(\mathbf x,\varepsilon)\cap {\rm int}(\texttt C(t_{j + 1})))\\
 & &
 \hspace{1.5cm}\geqslant
 \mathbb P(\|\texttt H_{j + 1} -\mathbf x\|\leqslant
 \varepsilon\wedge\varepsilon_{j + 1,\mathbf x})\\
 & &
 \hspace{1.5cm}\geqslant
 \mathbb P(\mathfrak c_1\|\overline{\tt X}_j -\mathbf x\|\leqslant 2^{-1}
 (\varepsilon\wedge\varepsilon_{j + 1,\mathbf x}),
 \mathfrak c_2\texttt R_{j + 1}(\mathbf x)\leqslant 2^{-1}
 (\varepsilon\wedge\varepsilon_{j + 1,\mathbf x}))\\
 & &
 \hspace{1.5cm}\geqslant
 \mathbb P(\|\overline{\tt X}_j -\mathbf x\|\leqslant (2\mathfrak c_1)^{-1}
 (\varepsilon\wedge\varepsilon_{j + 1,\mathbf x}),
 \texttt H_j\in {\rm int}(\texttt C(t_j)))\\
 & &
 \hspace{6.5cm}\times
 \mathbb P(\texttt R_{j + 1}(\mathbf x)\leqslant (2\mathfrak c_2)^{-1}
 (\varepsilon\wedge\varepsilon_{j + 1,\mathbf x})).
\end{eqnarray*}
On the one hand, since $\overline{\tt X}_j =\texttt H_j$ on $\{\texttt H_j\in {\rm int}(\texttt C(t_j))\}$, and since $\mathbf x\in\texttt C(t_{j + 1})\subset\texttt C(t_j)$, by (\ref{convergence_multidimensional_3}),
\begin{eqnarray*}
 & &
 \mathbb P(\|\overline{\tt X}_j -\mathbf x\|\leqslant (2\mathfrak c_1)^{-1}
 (\varepsilon\wedge\varepsilon_{j + 1,\mathbf x}),
 \texttt H_j\in {\rm int}(\texttt C(t_j)))\\
 & &
 \hspace{4cm} =
 \mathbb P_{\texttt H_j}(\overline B_d(\mathbf x,(2\mathfrak c_1)^{-1}
 (\varepsilon\wedge\varepsilon_{j + 1,\mathbf x}))\cap
 {\rm int}(\texttt C(t_j))) > 0.
\end{eqnarray*}
On the other hand, since $\sigma(\mathbf x)\in {\rm GL}_m(\mathbb R)$, and since $[\texttt Z_{j + 1}]_1,\dots,[\texttt Z_{j + 1}]_m$ are i.i.d. Gaussian random variables,
\begin{eqnarray*}
 & &
 \mathbb P(\texttt R_{j + 1}(\mathbf x)\leqslant (2\mathfrak c_2)^{-1}
 (\varepsilon\wedge\varepsilon_{j + 1,\mathbf x}))\\
 & &
 \hspace{1cm}\geqslant
 \mathbb P(\|\sigma(\mathbf x)\|_{\rm op}\|\texttt Z_{j + 1} +
 \sigma(\mathbf x)^{-1}b(\mathbf x)\Delta_n\|
 \leqslant (2\mathfrak c_2)^{-1}
 (\varepsilon\wedge\varepsilon_{j + 1,\mathbf x}))
 \geqslant\mu_{j + 1,\mathbf x}(\varepsilon)
\end{eqnarray*}
with
\begin{eqnarray*}
 \mu_{j + 1,\mathbf x}(\varepsilon)
 & = &
 \prod_{k = 1}^{m}
 \mathbb P(|[\texttt Z_{j + 1}]_k +
 [\sigma(\mathbf x)^{-1}b(\mathbf x)]_k\Delta_n|\\
 & &
 \hspace{3.5cm}
 \leqslant (2\mathfrak c_2\sqrt m\|\sigma(\mathbf x)\|_{\rm op})^{-1}
 (\varepsilon\wedge\varepsilon_{j + 1,\mathbf x})) > 0.
\end{eqnarray*}
{\bf Step 3 (conclusion).} By Step 2, for any $j\in J_n\backslash\{0\}$,
\begin{displaymath}
\mathbb P_{\texttt H_j}(\overline B_d(\mathbf x,\varepsilon)
\cap {\rm int}(\texttt C(t_j))) > 0\textrm{ $;$ }
\forall\varepsilon > 0\textrm{, }
\forall\mathbf x\in {\rm int}(\texttt C(t_j)).
\end{displaymath}
Thus, by Lemma \ref{lemma_convergence_multidimensional}, $\widehat{\tt C}_{N,j}$ converges pointwise to $\texttt C(t_j)$:
\begin{displaymath}
d(\mathbf x,\widehat{\tt C}_{N,j})
\xrightarrow[N\rightarrow\infty]{\mathbb P} 0
\textrm{ $;$ }\forall\mathbf x\in {\rm int}(\texttt C(t_j)).
\end{displaymath}
%


%
\subsubsection{Proof of Lemma \ref{lemma_convergence_multidimensional}}\label{section_proof_lemma_convergence_multidimensional}
Consider $\mathbf x\in C$ and
\begin{displaymath}
\delta_N(\mathbf x) :=\min_{i = 1}^{N}d(\mathbf x,X_i).
\end{displaymath}
Since $\{X_1,\dots,X_N\}\subset\widehat C_N$ by the definition of $\widehat C_N$,
\begin{eqnarray*}
 d(\mathbf x,\widehat C_N)
 & = &
 \inf_{\mathbf y\in\widehat C_N}d(\mathbf x,\mathbf y)\\
 & \leqslant &
 \min_{\mathbf y\in\{X_1,\dots,X_N\}}d(\mathbf x,\mathbf y) =
 \delta_N(\mathbf x).
\end{eqnarray*}
Moreover, for any $\varepsilon > 0$,
\begin{eqnarray*}
 \mathbb P(\delta_N(\mathbf x) >\varepsilon) & = &
 \mathbb P\left(\bigcap_{i = 1}^{N}\{d(\mathbf x,X_i) >\varepsilon\}\right)\\
 & = &
 \mathbb P(d(\mathbf x,X) >\varepsilon)^N =
 (1 - F_{\bf x}(\varepsilon))^N,
\end{eqnarray*}
where $F_{\bf x}$ is the distribution function of $d(\mathbf x,X)$, and
\begin{eqnarray*}
 F_{\bf x}(\varepsilon) & = &
 \mathbb P(d(\mathbf x,X)\leqslant\varepsilon,H\in {\rm int}(C)) +
 \mathbb P(d(\mathbf x,X)\leqslant\varepsilon,H\not\in {\rm int}(C))\\
 & \geqslant &
 \mathbb P(d(\mathbf x,H)\leqslant\varepsilon,H\in {\rm int}(C))\\
 & &
 \hspace{3cm}{\rm because}\quad
 \{H\in {\rm int}(C)\}\subset\{X = H\}.
\end{eqnarray*}
Thus,
\begin{displaymath}
\mathbb P(d(\mathbf x,\widehat C_N) >\varepsilon)\leqslant
(1 -\mathbb P_H(\overline B_d(\mathbf x,\varepsilon)\cap {\rm int}(C)))^N.
\end{displaymath}
In conclusion, by Inequality (\ref{lemma_convergence_multidimensional_1}),
\begin{displaymath}
d(\mathbf x,\widehat C_N)
\xrightarrow[N\rightarrow\infty]{\mathbb P} 0.
\end{displaymath}
%


%

%
\end{document}